\documentclass[english]{amsart}

\usepackage{babel}
\usepackage{amstext}
\usepackage{amsmath}
\usepackage{amsfonts}
\usepackage{latexsym}
\usepackage{ifthen}
\usepackage{xypic}
\xyoption{all}
\newcommand{\sO}{{\mathcal O}}
\newcommand{\I}{{\mathcal I}}
\newcommand{\PN}{{\mathbb P}}
\newcommand{\PF}{{\mathbb F}}
\newcommand{\KQ}{{\mathbb Q}}
\newcommand{\KZ}{{\mathbb Z}}

\newcommand{\FS}{{\mathbb F}}
\newcommand{\Pic}{{\rm Pic}}
\newcommand{\lra}{\longrightarrow}
\newcommand{\KC}{{\mathbb C}}

\newcommand{\Bl}{{\rm Bl}}
\newcommand{\exc}{{\rm exc}}
\newcommand{\X}{{\mathcal X}}
\newcommand{\F}{{\mathcal F}}
\newcommand{\E}{{\mathcal E}}
\newcommand{\Bs}{{\rm Bs}}

\newcommand{\bP}{{\mathbb P}}
\newcommand{\bQ}{{\mathbb Q}}
\newcommand{\sF}{{\mathcal F}}
\newcommand{\sI}{{\mathcal I}}
\newcommand{\sE}{{\mathcal E}}

\newcounter{lemma}

\newtheorem{lemma1}[lemma]{\setcounter{equation}{0}}

\newenvironment{lemma}{\begin{lemma1}{\bf Lemma.}}{\end{lemma1}}

\newenvironment{theorem}{\begin{lemma1}{\bf Theorem.}}{\end{lemma1}}

\newenvironment{theorem2}[1]{\begin{lemma1}{\bf Theorem [#1].}}{\end{lemma1}}
\newenvironment{proposition}{\begin{lemma1}{\bf
      Proposition.}}{\end{lemma1}}
\newenvironment{corollary}{\begin{lemma1}{\bf Corollary.}}{\end{lemma1}}
 
\newenvironment{remark}{\begin{lemma1}{\bf Remark.}\rm}{\end{lemma1}}

\begin{document}

\title {Almost del Pezzo manifolds}
\author[P. Jahnke]{Priska Jahnke} 
\address{P. Jahnke - 
Mathematisches Institut - Universit\"at Bayreuth - D-95440 Bayreuth, Germany} 
\email{priska.jahnke@uni-bayreuth.de}
\author[Th. Peternell]{Thomas Peternell}
\address{Th. Peternell - 
Mathematisches Institut - Universit\"at Bayreuth - D-95440 Bayreuth, Germany} 
\email{thomas.peternell@uni-bayreuth.de}
\thanks{Both authors gratefully acknowledge support by the Schwerpunkt program 
{\em Globale Methoden in der komplexen Geometrie} of the Deutsche Forschungsgemeinschaft.}
\date{\today}
\maketitle
\tableofcontents

\section{Introduction}
\setcounter{lemma}{0}

A del Pezzo manifold is a projective manifold $X$ of dimension $n$ whose anticanonical bundle is ample and divisible by $n-1$ in the 
Picard group. These manifolds are classical objects in algebraic geometry and completely classified (Iskovskikh, Fujita, ...). In terms of
differential geometry one classifies manifolds with positive Ricci
curvature whose canonical class has the above divisibility. It is therefore
natural to allow some degeneracies of the curvature and ask for a classification. This is the purpose of this paper: we consider projective 
manifolds $X$ with nef anticanonical class $-K_X$ such that $({-}K_X)^n > 0.$ In terms of differential geometry, the Ricci curvature is non-negative,
and the curvature is positive at some point. 
\vskip .2cm \noindent There is a relation to certain singular del Pezzo varieties: one can contract all curves in $X$ which are $K_X-$trivial and
obtains a birational map $\psi: X \to X'$, the ``anticanonical morphism'', with a singular Gorenstein variety $X'$ whose anticanonical bundle is ample 
and has the same divisibility. 
These singular del Pezzo varieties admit only a partial classification, on the other hand the existence of a {\it crepant} resolution 
(i.e. $K_X = \psi^*(K_{X'})$) is a priori a strong condition. 
\vskip .2cm \noindent Manifolds whose anticanonical bundles are big and nef are often called {\it almost Fano}, so we will speak in our context 
of {\it almost del Pezzo manifolds} - this explains the title. 
\vskip .2cm \noindent 
The paper is organized as follows. In section 3 we consider the case that $X$ has dimension 3 and Picard number 2. Here $X$ carries a unique 
Mori contraction which is either a quadric fibration, a $\bP_1-$bundle or the blow-up of a smooth point. Moreover we only need to treat the case that 
the anticanonical morphism is small, since the divisorial case was already treated in [JPR05]. The complete classification is given in the
theorems 3.1, 3.5, 3.6 and 3.7 below.
\vskip .2cm \noindent In case $\rho(X) \geq 3, $ but $X$ still of dimension 3, the paper [CJR06] plays an important role, and we show that either
after possibly performing a finite number of flops, $X$ is the blow-up of a certains number of points of a threefold classified in section 3, or 
$X$ is the projectivization of a rank 2-bundle over $\bP_2$, $\FS_2$ or $\bP_1 \times \bP_1,$ which can be written down explicitly (Theorem 4.1). 
\vskip .2cm \noindent In the last section we give the classification in dimension $n \geq 4,$ using the previous results. Here any Mori contraction
is either a $\bP_{n-2}-$bundle over a smooth surface, a quadric bundle over $\bP_1$ or the blow-up of a smooth point in another almost del Pezzo
manifold $Y.$ Using Fujita's partial classification of Gorenstein del Pezzo $n-$folds, we arrive at the classification given in the theorems 5.3, 5.8 and 5.9.

\section{Preliminaries}
\setcounter{lemma}{0}

Let $X$ be a smooth almost Fano manifold of dimension $n$ which is to say that $-K_X$ is big and nef. Suppose that $X$ is of index $n-1$, i.e. 
 \[-K_X = (n-1)H\] 
for some $H \in \Pic(X)$. Since $H$ is big and nef, the linear system $|mH|$ is base
point free for all $m \gg 0$. Since $|mH|$ and $|(m+1)H|$ define the
same map 
 \[\psi\colon X \lra X'\]
for $m \gg 0$, we find 
 \[H = \psi^*H'\]
for some $H' \in \Pic(X')$ with $-K_{X'} = (n-1)H'$, hence $X'$ is a singular del
Pezzo variety. We define the {\em degree} of $X$ (resp. of $X'$) to be
 \[d = H^n = (H')^n.\]
Smooth del Pezzo manifolds are classified by Fujita and Iskovskikh as follows.

\begin{theorem2}{\cite{Fu80}, \cite{Fu90}, \cite{Isk78}, \cite{Isk80}} \label{delpezzo}
Let $X$ be a del Pezzo manifold of dimension $n \geq 3.$ Then $X$ is one of the
following
 \begin{enumerate}
  \item $d = 1$, and $X \lra W$ is a double cover of the Veronese
    cone, ramified along a cubic,
  \item $d = 2$, and $X \lra \PN_n$ is a double cover, ramified along
    a quartic,
  \item $d = 3$, and $X \subset \PN_{n+1}$ is a cubic,
  \item $d = 4$, and $X \subset \PN_{n+2}$ is the complete intersection of
    two quadrics,
  \item $d = 5$, and $X \subset \PN_{n+3}$ is a linear section of the
    Grassmannian ${\rm Gr}(2,5) \subset \PN_9$ (embedded by
    Pl\"ucker). In particular $n \leq 6$,
  \item $d = 6$, and either
    \begin{enumerate}
     \item $X = \PN(T_{\PN_2})$,
     \item $X = \PN_1 \times \PN_1 \times \PN_1$ or $X = \bP_2 \times \bP_2$,
    \end{enumerate}
  \item $d = 7$, and $X = \Bl_p(\PN_3)$,
  \item $d = 8$, and $X = \PN_3$ with $H = \sO(2)$.  
 \end{enumerate}
\end{theorem2}

We now go back to the case that $X'$ is singular and Fano. Following the
notation in \cite{Fu90}, then $(X', H')$ is a polarized variety and we
obtain $\Delta(X', H') = 1$. By \cite{Shin} for $\dim(X') = 3$ and
\cite{Fu90} in general, $-K_{X'}$ is spanned, we may hence
assume that $\psi$ is the anticanonical map, i.e. the Stein
factorisation of the morphism defined by $|{-}K_X|$. Again by
\cite{Shin} and \cite{Fu90}, $H'$ is spanned for $(H')^n \ge
2$. For $(H')^n = 1$, the base locus of $|H'|$ is one point, contained in
$X'_{reg}$, hence $\Bs|H| = \Bs|H'|$ in any case. 

By the Riemann--Roch theorem and $H^i(X, H) = 0$ for $i > 0$ we get
 \[h^0(X, H) = H^n +n-1 = d+n-1.\]

\section{Threefolds with Picard number two}
\setcounter{lemma}{0}

Throughout this section we assume $\dim(X) = 3$, $\rho(X) = 2$ and $X$ not Fano, that is the first
non--trivial case. By Mori's classification (\cite{Mori}), any elementary extremal contraction 
 \[\phi\colon X \lra Y\] 
is either a del Pezzo fibration with general fiber $\PN_1 \times \PN_1$, a $\PN_1$--
bundle over a smooth surface, or the blowup of a smooth threefold in a
point. A more detailed description of possible contractions in our situation can be
found in \cite{CJR}. 

Since $K_X$ is not nef, by the cone theorem there
exists exactly one elementary extremal contraction $\phi$ and we
obtain a diagram
 \[\xymatrix{X \ar[r]^{\phi} \ar[d]^{\psi} & Y\\
             X'}\]
where now $X'$ is a Gorenstein Fano threefold with at most canonical
singularities and $\rho(X') = 1$. Two different cases may occur:
 \begin{enumerate}
  \item {\em The divisorial case:} $\psi$ contracts an irreducible
    divisor to a curve or a point. Then $X'$ is $\KQ$--factorial, but
    has canonical, non--terminal singularities.
  \item {\em The small case:} $\psi$ contracts only finitely many
    curves to points. Then $X'$ has terminal, but
    non--$\KQ$--factorial singularities.
\end{enumerate}

The divisorial case was already treated in the paper \cite{JPR}, the
result is 

\begin{theorem2}{\cite{JPR}}\label{div}
 Let $X$ be a smooth almost Fano threefold of index $2$ with $\rho(X)
 = 2$, such that $\psi$ is divisorial. Then $X$ is one of the
 following, and all of these cases really exist. The number in
 brackets refers to \cite{JPR}.
  \begin{enumerate}
    \item $X \to \PN_1$ is a del Pezzo fibration with general fiber
      $\PN_1 \times \PN_1$ and either
      \begin{enumerate} 
        \item $d = 1$, $X' \to W$ is a double cover of the Veronese cone, singular along a rational curve of degree $4$
      (A.2.12),
        \item $d = 2$, $X' \to \PN_3$ is a double cover, singular along a conic (A.2.15),
        \item $d = 2$, $X' \to \PN_3$ is a double cover, singular along an elliptic curve of degree $4$
      (A.2.9),
      \item $d = 4$, $X'\subset \PN(1^2, 2^3)$ is a hypersurface of degree $2$, singular along a conic (A.2.14).
    \end{enumerate}
    \item $X = \PN(\F)$ for some rank two vector bundle on $\PN_2$ and either
      \begin{enumerate}
        \item $d=3$, $X' \subset \PN_4$ is a cubic,
         singular along a line or the rational normal curve of degree $4$ and
         $\F \in {\mathcal M}(-1,4)$ is a stable Hulsbergen
          bundle (A.3.3 and A.3.4), 
         \item $d=6$, $X' \subset \PN_7$ is singular along a line and
           $\F$ is determined by $0 \lra \sO \lra \F \lra \I_p(-1)
           \lra 0$ (A.3.2),
         \item $d=9$, $X' = \PN(1^3,3)$ and $\F = \sO_{\PN_2} \oplus
           \sO_{\PN_2}(3)$ (A.3.1).
   \end{enumerate}
    \item $X = \Bl_p(V_{2,d+1})$ is the blowup in a
      point of a smooth del Pezzo threefold of degree $d+1$ and either
     \begin{enumerate}
       \item $d = 1$, $X' \to W$ is a double
         cover of the Veronese cone, singular along a conic or a
         smooth curve of degree $8$ and genus $3$ (A.5.5 and
         A.5.6),
         \item $d = 2$, $X' \subset \PN_5$ is singular along an
           elliptic curve of degree $6$ (A.5.7).     
     \end{enumerate}
  \end{enumerate}
\end{theorem2}

\

From now on we assume {\bf $\psi$ is small.} Then by \cite{Kollar}, there exists the following flop--diagram 
\[\xymatrix{X \ar@{-->}[rr]^{\chi} \ar[dr]^{\psi} \ar[d]_{\phi} & & X^+
      \ar[dl]_{\psi^+} \ar[d]^{\phi^+}\\
       Y & X' & Y^+}\]
where the rational map $\chi$ is an isomorphism outside the
exceptional locus of $\psi$ and $X^+$ is again a smooth almost
Fano threefold with anticanonical map $\psi^+$ and extremal
contraction $\phi^+$. Our assumption $\rho(X) = 2$ implies that $\chi$ does
not depend on the choice of some $\psi$--negative divisor in $X$. The
index of $X^+$ is again $2$, i.e. $-K_{X^+} = 2H^+$ for some $H^+ \in \Pic(X^+)$.

\vspace{0.2cm}

The following Lemma is essentially \cite{AG5}, Remark~4.1.10:

\begin{lemma}\label{hyp}
  If $d \le 2$, then $X \simeq X^+$ as abstract varieties.
\end{lemma}

\begin{proof}
First note that $X'$ is a double cover of some $\KQ$--factorial
threefold $W$: if $(H')^3 = 1$, then $|{-}K_{X'}|$ defines a double cover
of the Veronese cone, if $(H')^3 = 2$, then $X'$ is a double cover of
$\PN_3$, defined by $|H'|$.

Denote the birational involution induced on $X$ by $\sigma$ and let
$D$ be some divisor on $X$. Denote the strict transform under $\sigma$
by $D^{\sigma}$. Then $D + D^{\sigma}$ is the pullback of some
$\sigma$--invariant (Weil-) divisor $B'$ on $X'$ which actually comes from
$W$. As $W$ is $\KQ$--factorial, $mB'$ is Cartier. Then
  \[(D + D^{\sigma})\cdot l_{\psi} = \frac{1}{m} \psi^*(mB')\cdot l_{\psi} = 0\]
for any curve $l_\psi$ contracted by $\psi$. But then $D\cdot l_{\psi} =
-D^{\sigma}\cdot l_{\psi}$. This implies $\sigma: X \dasharrow X$ is the
flop (\cite{Kollar}) and in particular $X^+ \simeq X$.
\end{proof}

\vspace{0.2cm}

A {\em smoothing} of a singular Fano threefold $X'$ is a flat family
 \[\X \lra \Delta\]
over the unit disc, such that $\X_0 \simeq X'$ and $\X_t$ is a smooth
Fano threefold for $t \not= 0$. Namikawa has shown in \cite{Namikawa}
that a smoothing always exists if $X'$ has only terminal Gorenstein
singularities, not necessarily $\KQ$--factorial: in this case the
Picard groups of $X'$ and the general $\X_t$ are isomorphic (over
$\KZ$) by \cite{smoothing}. 

\begin{theorem2}{\cite{Namikawa}, \cite{smoothing}} \label{smoothing}
Let $X'$ be a Gorenstein Fano threefold with only
  terminal singularities (not necessarily $\bQ-$factorial). Then $X' $
  has a smoothing $\X \to \Delta$ and $\Pic(X') \simeq \Pic(\X_t)$. In particular, $X'$ and $\X_t$ have the
same Picard number, the same index and the same degree.
\end{theorem2}

\begin{corollary} \label{deg}
 If $\psi$ is small and $\rho(X) = 2$, then $1 \le H^3 \le 5$.
\end{corollary}

\begin{proof}
By Theorem~\ref{smoothing} above, $(\X_t, H_t)$ is a smooth del Pezzo
threefold of Picard number one, hence either $\X_t = \PN_3$ or $1 \le
H_t^3 \le 5$. But $\X_t = \PN_3$ implies the index of $X'$ is
$4$. Then $X' \simeq \PN_3$ by \cite{Shin}, which is impossible.
\end{proof}

\vspace{0.2cm}

The aim now is to describe all possible tripels $(X, X^+, \X_t)$ in
terms of their Mori fiber space structure. We will consider all
possibilities for $\phi$ seperately.

\

\subsection*{Case A: Del Pezzo fibrations.} Assume first $\phi\colon X \to
\PN_1$ is a del Pezzo fibration. Since $X$ has index $2,$ so does the general fiber $F$, hence $F \simeq \PN_1 \times \PN_1$ and actually
every fiber is a smooth quadric or a quadric cone. Define
  \[\E = \phi_*(H) = \oplus_{k=1}^4 \sO_{\PN_1}(a_k), \quad a_1 \le
  \dots \le a_4;\]
$\E$ is a vector bundle on $\PN_1$ of rank $r = h^0(F, H|_F) =
h^0(\PN_1 \times \PN_1, \sO(1,1)) = 4$. 
From $H^1(X, H) = 0$ we get $H^1(\PN_1, \E) = 0$. Hence
\[a_1 \ge -1.\]
By a result of Andreatta-Ballico-Wisniewski [BS95,p.338], the canonical map
 \[\phi^*\E \lra \sO_X(H)\]
is an epimorphism and therefore yields an embedding
$$ X \subset \bP(\sE) $$ such that $H = \sO_{\bP(\sE)}(1) \vert X.$ In other words, $X$ is a conic bundle. 

\begin{theorem} \label{dP}
Assume $\rho(X) = 2$, $\psi$ is small and $\phi$ is a del Pezzo
fibration. Then $X$ is a quadric bundle and belongs to the following list. 
 \begin{enumerate}
   \item $X \subset \PN_3 \times \PN_1$ from $|(2, 2)|$, here $d =
      2$, $X^+ \simeq X$ and $\X_t \to \PN_3$ is a double cover, 
   \item $X \subset \FS(0^3, 1)$ from $|\sO(2)+F|$, here $d = 3$,
   $X^+ = \Bl_p(V_{2,4})$ and $\X_t \simeq V_{2,3}$ (this is (3) in Theorem~\ref{bir}),
   \item $X \subset \FS(0^2, 1^2)$ from $|\sO(2)|$, here $d = 4$,
     $X^+$ is of the same type and $\X_t \simeq V_{2,4}$,
   \item $X \subset \FS(0, 1^3)$ from $|\sO(2)-F|$, here $d = 5$, $X^+ = \PN(\F)$ with some stable rank two bundle $\F \in
     {\mathcal M}(-1,2)$ (this is (1) in Theorem~\ref{cb}), and
     $\X_t \simeq V_{2,5}$,
    \item $X \subset \FS(-1, 0^2, 1)$ from $|\sO(2)+2F|$, here $d =
      2$, $X^+ \simeq X$ and $\X_t \to \PN_3$ is a double cover,
    \item $X \subset \FS(-1, 0^3)$ from $|\sO(2)+3F|$, here $d = 1$,
      $X^+ \simeq X$ and $\X_t \to W$ is a double cover of the
      Veronese cone.
  \end{enumerate}
\end{theorem}

\begin{proof}
We consider the embedding $X \subset \bP(\sE) = \FS$ and denote the tautological line bundle by $\sO_{\FS}(1)$. For some $\alpha \in \KZ$
  \[X \in |\sO_{\FS}(2) + \alpha F|.\]
As $H = \sO_{\FS}(1)|_X$ we obtain from the adjunction formula
  \begin{equation} \label{aalpha} 
     a_1 + a_2 + a_3 + a_4 - 2 + \alpha = 0
  \end{equation}
and  we find
  \[1 \le H^3 = 2(a_1 + a_2 + a_3 + a_4) + \alpha \le 5\]
by Corollary~\ref{deg}. Putting things together gives
 \begin{equation} \label{sumai}
    -1 \le a_1 + a_2 + a_3 + a_4 \le 3.
 \end{equation}
Then $a_1 = -1$ or $a_1 = 0$. 
\vskip .2cm \noindent {\bf (A)} First assume $\Bs|\sO_{\FS}(1)| = \emptyset$, i.e., $a_1
= 0$. We may assume that $a_4 > 0, $ hence $a_4 = 1;$ otherwise we are
clearly in case (1) of the theorem. \\
Now suppose furthermore that $a_3 = 0$.
Then we find a unique section $D \in |\sO_{\FS}(1) - a_4F|$ contracted to
$\PN_2$ by $|\sO_{\FS}(1)|$. As $\psi$ is not divisorial,
  \[\sO_{\FS}(1)^2\cdot D\cdot X = 2 - a_4 > 0,\]
hence $a_4 = 1$. Writing $s, t$ for homogeneous coordinates on $\PN_1$ and 
  \[s^2Q_1 + stQ_2+t^2Q_3 = 0\]
for the equation defining $X$, where $Q_i$ are general quadrics on $\PN_3$, we
see that $\psi$ contracts the $\PN_1$'s lying over the eight points of
intersection $Q_1 \cap Q_2 \cap Q_3$.
The map associated with $|H|$ is given by projection onto $\PN_3$ and $\psi$
is small if $X$ contains horizontal $\PN_1$'s. Notice that $\alpha  = 2 - a_4 = 1$. Here $D \simeq
\PN_2 \times \PN_1$ and for a general choice, $X \cap D \in |(2, 1)|$ is
a generic section over $\PN_2$,
isomorphic to $\PN_2$ blown up in four points. The corresponding
$(-1)$--curves are contracted by $\psi$, this is case (2) of the theorem. \\
In the case $a_3 > 0$, we get from \eqref{sumai}
  \[(a_1, a_2, a_3, a_4)_{\alpha} \in \{(0, 0, 1, 1)_0, (0, 0, 1, 2)_{-1},
  (0, 1, 1, 1)_{-1}\}.\]
In the first case we realize (3) in the theorem, in the second case, $\psi$ is divisorial since $X \in \vert  \sO_{\FS}(2) - F \vert  $ contains the exceptional locus
$\bP(\sO^2)$ of the map associated with $\vert \sO_{\FS}(1) \vert.$ 
In the last case $\alpha =
-1$ guarantees that the trivial
section contracted by $\sO_{\FS}(1)$ is in $X$. This leads to (4). 

\vspace{0.2cm} \noindent
{\bf (B)} Now assume $\Bs|\sO_{\FS}(1)| \not= \emptyset$, so that $a_1 =  -1.$ \\
The map 
$$H^0(\FS,
\sO_{\FS}(1)) \lra H^0(X, H)$$ is surjective,
since
$$H^1(\FS, \sO_{\FS}(1) - X) = 0.$$
Therefore $\Bs|\sO_{\FS}(1)| \cap X$
is at most a single point. Then $a_1 = -1$, $a_2 \ge 0$. Let $l =
\Bs|\sO_{\FS}(1)|$ be the corresponding section of $\FS$. Then $0 \le
X \cdot l \le 1$ and hence $0 \le -2 + \alpha \le 1$. By \eqref{aalpha}
  \[(a_1, a_2, a_3, a_4)_\alpha \in \{(-1, 0, 0, 1)_2, (-1, 0, 0, 0)_3\}.\]

\vskip .2cm \noindent 
{\bf (C) Existence and determination of the flop.} \\
It remains to show the existence of all cases and determine the
corresponding type of $X^+$. The type of the smoothing $\X_t$ follows
by Iskovskikh's list.

\vspace{0.2cm}

\noindent (1) In the case $X \subset \PN_1 \times \PN_3$ from $|\sO(2, 2)|$ the map
  corresponding to $|H|$ is the projection onto $\PN_3$ implying $X
  \simeq X^+$ by Lemma~\ref{hyp}. The involution $\sigma$ is given by
    \[[s:t] \mapsto [tQ_3 : sQ_1]\]
   times the identity.

\vspace{0.2cm}

\noindent (2) In the case $X \subset \FS(0^3, 1) = \Bl_{\PN_2}(\PN_4)$ we have
$X' \subset \PN_4$ a cubic given by
  \[x_0q_0 + x_1q_1\]
where $q_0, q_2$ are two general quadrics and $x_0, \dots, x_4$ are
homogeneous coordinates of $\PN_4$. The two quadrics intersect
the plane $\PN_2$ given by $I(x_0, x_1)$ in four points in $X'_{sing}$. 
The blowup of $q_0 \cap q_1$ gives $X$ in the same way, so this is
not the flop.

Our threefold $X$ intersects the exceptional divisor $\simeq
\PN_2 \times \PN_1$ of $\FS(0^3, 1) = \Bl_{\PN_2}(\PN_4)$ in a surface $E^+ \simeq \Bl_{p_1, \dots,
  p_4}(\PN_2) \in |(2, 1)|$. Restricted to $E^+$ the system $|H|$ corresponds to the
pullback of $\sO_{\PN_2}(1)$. The map corresponding to $|H|$ therefore 
contracts the four $(-1)$--curves in $E^+$. Then
  \[N_{l_{\psi}/X} = \sO_{\PN_1}(-1) \oplus  \sO_{\PN_1}(-1).\]
Let $\FS = \FS(0^3, 1)$. The exact sequence
  \[0 \lra \sO_{\FS}(-2F) \lra \sO_{\FS}(2H-F) \lra \sO_X(2H-F) \lra 0\]
shows $h^0(\FS, 2H-F) = 5$ but $h^0(X, 2H-F) = 6$. We have $\Bs|2H-F| = \exc(\psi)$. Let
$\tilde{X} = \Bl_{\exc(\psi)}(X)$. Denote the four exceptional divisors
$\simeq \PN_1 \times \PN_1$ by $E_i$. The system $|2H-F-\sum E_i|$
is spanned on $\tilde{X}$. It contracts the $E_i$'s in the other
direction. We get an induced map
  \[X^+ \lra \PN_5.\]
The difference $h^0(\FS, 2H-F) = 5$ but $h^0(X, 2H-F) = 6$ has the
following meaning. The five sections coming from $\FS$ have the form
  \[wx, wy, wz, sw^2, tw^2,\]
where $x,y,z$ correspond to the three trivial summands $\sO_{\PN_1}$
and where $w$ corresponds to $\sO_{\PN_1}(1)$. In other words they are
all reducible of the form $H + E$. The existence of one additional
section means: if we project $\PN_5$ onto $\PN_4$ from $[0:0:0:0:1]$
then we get $X' \subset \PN_4$ (divide by $w$).

Also a direct computation shows $(2H-F-\sum E_i)^3 = 4$. Then $X^+$ is
mapped onto some threefold $Y^+ \subset \PN_5$. On the strict
transform of $E^+$ in $\tilde{X}$ we find (using $K_{E^+} = -3H + \sum
E_i$ and $K_{E^+} + H = (-1, -1) + (1, 0) = -F$):
  \[2H-F-\sum E_i|_{E^+} = 2H + K_{E^+} + H - \sum E_i|_{E^+} = 0.\]
Then $E^+$ is contracted to a point. The image of $E^+$ in $X^+$ is
isomorphic to $\PN_2$ and mapped to a point.

\vspace{0.2cm}

\noindent (3) In the case $X \subset \FS(0^2, 1^2)$ from $|2H|$ we first note that $\FS(0^2,
1^2)$ is a small resolution of the double cone over the quadric $Q_2$.
 \[\xymatrix{ & \hat{\FS} = \PN(\sO_{\PN_1 \times \PN_1}^{\oplus 2} \oplus \sO_{\PN_1 \times
     \PN_1}(1, 1)) \ar[ld]\ar[rd] & \\
      \FS(0^2, 1^2)\ar[rd] & & \FS(0^2, 1^2) \ar[ld]
      \\
    & \widehat{\widehat{\PN_1 \times \PN_1}} &}\]
Then $X' = Q_4 \cap \widehat{\widehat{Q}}_2$ with $Q_4$ general. 

\vspace{0.2cm}

\noindent (4) The case $X \subset \FS(0, 1^3)$ from $|2H - F|$ where $H =
  \sO_{\FS}(1)$. Here we have one $-K_X = 2H$--trivial curve $l_{\psi}$
  cut out by the three sections of $|H - F|$. Note that $|H|$ maps
  $\FS$ onto the cone over $\PN_1 \times \PN_2$ embedded by the Segre
  embedding. Blowing up the vertex of the cone we obtain a divisorial
  resolution $\hat{\FS} = \PN(\sO_{\PN_1 \times \PN_2} \oplus \sO_{\PN_1 \times
     \PN_2}(1, 1))$ with exceptional divisor $E \simeq \PN_1 \times
   \PN_2$. Blowing down $E$ in either direction first, we get two
   small resolutions:  
   \[\xymatrix{ & \hat{\FS} = \PN(\sO_{\PN_1 \times \PN_2} \oplus \sO_{\PN_1 \times
     \PN_2}(1, 1)) \ar[ld]\ar[rd] & \\
      \FS = \FS(0, 1^3)\ar[rd] & & \PN(\sO_{\PN_2} \oplus
      \sO_{\PN_2}(1)^{\oplus 3}) \ar[ld]
      \\
    & \widehat{\PN_1 \times \PN_2} &}\]
  The pullback of the tautological systems $H$ and $H^+$ respectively give the tautological
  system on $\hat{\FS}$ which
  we denote by $\hat{H}$. We have $\hat{X} \in |\hat{H} + (0,
  1)|$ and hence $X^+ \in |H^+ + F^+|$ where $F^+ \simeq \PN_2$ is a
  fiber. In particular $-K_{\hat{X}} = (\hat{H} + (1,1))_{\hat{X}}$ and
  $\hat{X}$ is Fano. 

Since $X^+|_{F^+}$ is a line, the induced projection map $X^+
  \to \PN_2$ is a $\PN_1$--bundle, i.e. $X^+ = \PN(\F)$ for some rank
  two bundle $\F$ on $\PN_2$. Assuming $\F$ to be normalized we
  compute $c_1(\F) = -1$ and $c_2(\F) = 2$. This is case (1) in the
  following Theorem~\ref{cb}.

\vspace{0.2cm}

\noindent (5) Here $|H'|$ is basepoint free, defining a double cover $X' \to
\PN_3$. Hence $X^+ \simeq X$ by Lemma~\ref{hyp}.

\vspace{0.2cm}

\noindent (6) In the last case $X$ is hyperelliptic, meaning $-K_{X'} = 2H'$ is
generated, defining a double cover of $X'$ onto the Veronese cone. Hence $X^+ \simeq X$ by Lemma~\ref{hyp}.
\end{proof}

\

\subsection*{Case B: Conic bundles} 
Assume now $\phi$ is a conic
bundle $X \to \PN_2$ with discriminant $\Delta$. Since $-K_X$ is divisible in $\Pic(X)$, there cannot be any
reducible fibers, hence $\Delta = \emptyset$ and $\phi$ is a $\PN_1$--bundle. 

\begin{theorem} \label{cb}
Assume $\rho(X) = 2$, $\psi$ is small and $\phi$ is a
$\PN_1$--bundle. Then $X = \PN(\F)$ with a stable rank $2$ bundle on
$\PN_2$ with $c_1(\F) = -1$ and $2 \le c_2(\F) \le 5$. Moreover,
$\F(2)$ is nef, but not ample and has only finitely many jumping
lines. We have
 \begin{enumerate}
   \item $c_2(\F) = 2$. Then $d = 5$, $X^+$ admits a del Pezzo
     fibration as in (4) of Theorem~\ref{dP} and $\X_t \simeq V_{2,5}$,
   \item $c_2(\F) = 3$. Then $d = 4$, $X^+ = \Bl_p(V_{2,5})$ and $\X_t \simeq V_{2,4}$,
   \item $c_2(\F) = 4$. Then $d = 3$, $X^+$ is of the same type, and $\X_t \simeq V_{2,3}$,
   \item $c_2(\F) = 5$. Then $d = 2$, $X^+ \simeq X$, and $\X_t \to
     \PN_3$ is a double cover. 
  \end{enumerate}
\end{theorem}

\begin{proof} We write $X = \bP(\sF)$ with $\eta = \sO_{\bP(\sF)}(1)$
  and normalize $\sF$ such that $c_1(\sF) = 0,-1.$ Let $L =
  \phi^*\sO_{\PN_2}(1)$. Then $-K_X = 2\eta + (3-c_1(\sF))L$ divisible
  implies $c_1(\sF) = -1$ and 
 \[-K_X = 2\eta + 4L.\]
It follows that $\sF(2)$ is nef but not ample. It is clear that 
$$ \sF \vert l = \sO \oplus \sO(-1)$$
for all but finitely many lines in $\bP_2.$ 
If this would hold for all lines, then $\sF = \sO \oplus \sO(-1)$ or $\sF = T_{\bP_2}(-2),$ see e.g. \cite{OSS}. But in both cases $\sF(2)$
would be ample, hence we must have splitting lines $l$ with
$$ \sF \vert l =  \sO(1) \oplus \sO(-2).$$ 
Now $X'$ is a del Pezzo variety with $\rho(X') = 1$ and so do the
smoothings $\X'_t.$ From the classification of smooth del Pezzo
threefolds and the cohomology of the ample generator we obtain 
$$ 3 \leq h^0( \eta + 2L) \leq 7,$$ 
so that
$$ 3 \leq h^0(\sF(2)) \leq 7.$$
Observe also that $H^q(\sF(2)) = 0$ for $q \geq 1$ so that Riemann-Roch yields
$ 3 \leq 9 - c_2(\sF) \leq 7,$ so that
$$ 2 \leq c_2(\sF) \leq 5.$$
Finally stability of $\sF$ is obvious since $h^0(\sF) = 0.$ 

\

We need to show the existence of the four cases, and to describe the
flops. Concerning existence let $\sF$ be a general member of the moduli space of stable rank 2-bundles on $\bP_2$ with $c_1(\sF) = -1$ and $2 \leq c_2(\sF) \leq 5.$ 
By \cite{LeP79}, $\sF(2)$ is spanned. Moreover $\sF$ has only finitely
many jumping lines by \cite{Hulek}. Set $X = \bP(\sF)$. Then $-K_X$ is
spanned and big and the map associated with $|{-}K_X|$ is
small. We denote $F = \phi^*\sO_{\PN_2}(1)$.

\vspace{0.2cm}

\noindent (1) Assume $c_2 = 2$. By \cite{Hulek}, Proposition~8.2 there is
exactly one jumping line. Note $h^0(\F(2)) = d+2 = 7$. Tensorising the ideal sequence of a general line in
$\PN_2$ with $\F(2)$ then shows $h^0(\F(1)) \ge 2$, i.e. $\F$ is a Hulsbergen
bundle. Then $\F$ is determined by an extension
 \begin{equation} \label{huls}
   0 \lra \sO \lra \F(1) \lra \I_Y(1) \lra 0, 
\end{equation}
where Y consists of two general points. The jumping line is the unique
line through these points and sequence \eqref{huls} implies $h^0(\F(1))
= 2$. The base locus of $|H-F|$ is exactly the exceptional locus of
$\psi$, which is the minimal section $C_0$ over the jumping
line. Flopping $C_0$, the system becomes base point free, hence $X^+$
admits a del Pezzo fibration. This must be case (4) in
Theorem~\ref{dP} for numerical reasons.

\vspace{0.2cm}

\noindent (2) Assume $c_2 = 3$. As in (1) we find $h^0(\F(1)) \ge 1$, hence $\F$
is determined by an extension \eqref{huls}, where now $Y = \{p_1, p_2,
p_3\}$. Since $\F(2)$ is supposed to be nef, the points are in general
position, hence $\F$ is a Hulsbergen bundle and $h^0(\F(1)) = 1$. We
have exactly $3$ jumping lines, each connecting $2$ of the $3$
points. The unique section $S \in |H-F|$ is a smooth del Pezzo
surface of degree $6$, containing $\exc(\psi) = \{C_1, C_2, C_3\}$ as
$({-}1)$--curves. Flopping the curves in $X$ means blowing them down
in $S$, hence the image $S^+$ becomes a contractible $\PN_2$. This
shows $X^+ = \Bl_p(V_{2,5})$.

\vspace{0.2cm}

\noindent (3) Assume $c_2 = 4$. If $h^0(\F(1)) \not= 0$, then $|H|$ contracts
the unique section of $|H-F|$ to a curve, meaning $\psi$ is
divisorial. Hence $h^0(\F(1)) = 0$. We claim $Y^+ = \PN(\F^+)$ with
$\F^+ \in {\mathcal M}(-1,4)$ is of the same type. The twisted ideal sequence of a line in $\PN_2$ gives $h^0((S^2\F)(3))
\ge 3$ and 
 \[\Bs|2H-F| = \exc(\psi).\] 
So after flop $|2H-F|$ becomes free, it hence remains to show that the induced map
$\phi^+$ indeed maps $X^+$ onto $\PN_2$. 

Since $X' \subset \PN_4$ is a cubic, a general member $S \in |H'|$ is
a smooth cubic, not meeting $X'_{sing}$. We may hence identify $S$
with its pullbacks to $X$ and $X^+$. Then 
 \[\phi|_S\colon S \lra \PN_2\]
is the blowup of $6$ general points $p_1, \dots, p_6$. We find that
$(2H-F)|_S$ is spanned and big, contracting the $6$ conics,
each through $5$ of the points $p_1, \dots, p_6$. Hence 
 \[\phi^+|_{S^+}\colon S^+ \lra \PN_2\]
is again birational. Assume $\phi^+$ is birational. Then $X^+ =
\Bl_p(V_{2,4})$ by classification, but $V_{2,4}$ does not contain a
family of $\PN_2$'s. This shows $\phi^+\colon X^+ \to \PN_2$ as
claimed.

\vspace{0.2cm}

\noindent (4) If $c_2 = 5$, then $|H'|$ is base point free, defining $X' \to
\PN_3$ a double cover. Then $X^+ \simeq X$ by Lemma~\ref{hyp}.
\end{proof}

\

\subsection*{Case C: Blowups} If $\phi$ is birational, then
$-K_X = 2H$ implies $\phi$ is the blowup of a smooth del Pezzo threefold
$Y$ of Picard number $1$ in a point. By classification then $Y =
V_{2,d+1}$, a smooth del Pezzo threefold of degree $d+1$ as in
Iskovskikh's list.

\begin{theorem} \label{bir}
Assume $\rho(X) = 2$, $\psi$ is small and $\phi$ is birational. Then
$X = \Bl_p(Y)$ for a general point $p$ in a smooth del Pezzo threefold $Y = V_{2,d+1}$, such that
 \begin{enumerate}
   \item $d = 1$, $X^+ \simeq X$ and $\X_t
     \to W$ is a double cover of the Veronese cone,
   \item $d = 2$, $X^+ \simeq X$ and $\X_t \to \PN_3$
     is a double cover,
   \item $d = 3$, $X^+$ admits a del Pezzo fibration
     as in (2) of Theorem~\ref{dP}, and $\X_t \simeq V_{2,3}$,
   \item $d = 4$, $X^+ = \PN(\F)$ as in (2) of
     Theorem~\ref{cb}, and $\X_t \simeq V_{2,4}$. 
  \end{enumerate}
\end{theorem}

\begin{proof}
Consider $Y = V_{2,d+1}$ and let $p \in Y$ be a general point. Then
there are only finitely many lines through $p$. Their strict
transforms in $X$ are the exceptional curves of $\psi$. This shows the
existence of (1)--(4). Lemma~\ref{hyp} implies $X \simeq X^+$ in (1)
and (2); the structure of $X^+$ in (3) and (4) follows by
the Theorems~\ref{dP} and \ref{cb} above.
\end{proof}

\section{The general case in dimension $3$}
\setcounter{lemma}{0}

Let now $X$ be any smooth almost Fano threefold of index two which is not Fano. Here we do not make any restriction on the Picard number. The
following theorem uses in a very essential way Proposition~2.10 in
\cite{CJR} and describes $X$ up to flops:

\begin{theorem} \label{gen}
Let $X$ be a smooth almost Fano threefold of index $2$. Then either
  \begin{enumerate}
   \item There exists a finite sequence of flops $X \dasharrow
     \Bl_{p_1, \dots, p_r}(X_0)$, where $X_0$ is a smooth almost Fano
     threefold (possibly Fano) of index $2$ with $\rho(X_0) \le 2$. If we write $-K_{X_0} =
     2H_0$, then $r < H_0^3$,
   \item $X = \PN(\F)$ is a $\PN_1$--bundle over $S = \PF_2,$ or $ \PN_1 \times \PN_1$, and $\F$ is a nef
     rank $2$ vector bundle with $c_1(\F) = -K_S$ and $0 \le c_2(\F)
     \le 7$, $c_2(\F) \not=1$, given by an extension
   \[0 \lra \sO_S \lra \F \lra \sO_S(-K_S) \otimes \I_Z \lra 0,\]
   where $\I_Z$ is the ideal sheaf of $c_2(\F)$ points on $S$ which are almost in general position in the following sense. Two points of $Z$ are
on a ruling line of one of the two rulings in case $S = \bP_1 \times \bP_1$ resp. on a ruling line of the unique ruling in case $S = \PF_2.$ The other points 
are in general position. 
 \end{enumerate}
All these cases really exist.
\end{theorem}

\begin{remark}
  \begin{enumerate}
    \item The description in (1) is in general not unique. For example
      $\Bl_p(\PN_3) = \PN(\sO_{\PN_2} \oplus \sO_{\PN_2}(1))$.
    \item The number of flops in the description (1) is by
      construction bounded by the number of blowups $r$. Flops only
      occur if $X_0 = \PN(\F)$.
   \item By \cite{CJR}, the Picard number of $X$ is effectively
     bounded by $10$; the Picard number of an anticanonical model
     $X'$ is bounded by $3$, with equality only for $X = X' = \PN_1
     \times \PN_1 \times \PN_1$.   
  \end{enumerate}
\end{remark}

\begin{proof}[Proof of Theorem~\ref{gen}.] Let $X$ be a smooth almost
  Fano threefold of index $2$. Then $-K_X\cdot C \ge 2$ for any rational
  curve of positive degree, i.e. the pseudo--index $i_X > 1$. Then
  Proposition~2.10 in \cite{CJR} applies: we are either in case (1) of the theorem, or $\rho(X) =
  3$ and $X = \PN(\F)$ over $\PF_2$ or $\PN_1 \times \PN_1$. It
  remains hence to describe the bundle $\F$ in the second case.
We present two completely different methods, one is birational and reduces the problem to the base space
$\bP_2;$ the other is more vector-bundle-theoretic. 

\subsection*{The Birational Method}
Almost Fano $\PN_1$--bundles over $\PN_2$ are classified in the last section and
 Theorem~\ref{delpezzo}. The aim is hence to reduce the problem from the general $S$ to
 $\PN_2$ by blowing up and down.

\vspace{0.2cm}

We generalise the situation and consider rank two vector bundles $\F$
 on smooth almost Fano surfaces $S$, fitting into a sequence
  \begin{equation} \label{ext0}
   0 \lra \sO \lra \F \lra \sO(-K_S) \otimes \I_{\F} \lra 0,
  \end{equation}
where $\I_{\F}$ is the ideal sheaf of $c_2(\F)$ points, not
necessarily in general position. We claim

\begin{lemma}
 1.) Let $\F_0$ be a vector bundle of type \eqref{ext0} on the surface
 $S_0$. Let $q$ be a general point on $X_0 = \PN(\F_0)$. Then there exists a diagram
 \begin{equation} \label{blowup}
   \xymatrix{X_1 = \PN(\F_1) \ar@{<-->}[r]^{\hspace{0.5cm} flop} \ar[d]^{f_1} & X^+_0
   \ar[r]^{\hspace{-0.5cm} \Bl_q} & \PN(\F_0) = X_0 \ar[d]^{f_0}\\
            S_1 \ar[rr]^{\pi=\Bl_q} && S_0}
 \end{equation}
such that $\F_1$ is again of type \eqref{ext0} and $\I_{\F_1}$
corresponds to the same points as $\I_{\F_0}$. 

\vspace{0.2cm}

2.) Conversely, let $\F_1$ be a vector bundle of type \eqref{ext0} on the surface
 $S_1$ and $X_1 = \PN(\F_1)$. Assume $S_1$ is not minimal and let $\pi\colon S_1 \to S_0$ be the
 contraction of a $({-}1)$--curve to a point $q$, not contained in $Z$. Then there exists the diagram
 \eqref{blowup}, where $\F_0$ is of type \eqref{ext0} and $\I_{\F_0}$
corresponds to the same points as $\I_{\F_1}$.
\end{lemma}

In other words: the first Chern classes of $\F_i$ both equal the
anticanonical divisor of the underlying surface $S_i$, and the second
Chern classes remain unchanged under blowup. 

\begin{proof}
1.) Let $Z_0$ be the support of $\I_{\F_0}$. We denote the image of $q$ in $S_0$ again by $q$. Then $q$ general
implies $q \not\in Z_0$. The general construction is now well known: the strict transform $C_0 \subset
X^+_0$ of the fiber in $X_0$ containing $q$ is an anticanonically trivial
curve with normal bundle of type $({-}1,{-}1)$. Blowing up $X^+_0$ along
$C_0$ and then again down in the other direction gives the flop to
$X_1$. Denote the the image of the flopping curve by $C_1$.

Denote by $E_0 \simeq \PN_2$ the exceptional divisor of $X^+_0 \to
X_0$. Then $E_0$ meets $C_0$ transversally in a single point and the strict transform $E_1$ of $E_0$ in $X_1$ is isomorphic to
$\PF_1$, containing $C_1$ as its minimal section. The image $f_1(E_1)$ in
$S_1$ is exactly the exceptional curve of the blowup $\pi$.

This proves $X_1 = \PN(\F_1)$ for some rank two vector bundle $\F_1$, we
have to show the existence of the sequence \eqref{ext0} for $\F_1$
with $\I_{\F_1}$ as claimed. To this end we chase a general section $H_0 \in |\sO_{X_0}(1)|$ through the
diagram. Then sequence \eqref{ext0} for $\F_0$ reads  
 \[H_0 = \Bl_{Z_0}(S_0), \quad N_{H_0/X_0} = -K_{H_0}.\] 
Let $H^+ \simeq H_0$ be the strict transform of $H_0$ in $X^+_0$ and $H_1
\in |\sO_{X_1}(1)|$ its strict transform in $X_1$. Since $H^+$ meets the flopping curve $C_0$
transversally in one point $p$, we get $H_1 = \Bl_p(H^+)$ with
exceptional curve $C_1$. Since $C_1$ is a section over the exceptional
curve of $\pi\colon S_1 \to S_0$, we find 
 \[H_1 = \Bl_{Z_1}(S_1), \quad N_{H_1/X_1} = -K_{H_1},\] 
where $Z_1 = \pi^{-1}(Z_0) \simeq Z_0$. This shows $\F_1$ is of type
\eqref{ext0} and proves 1.).

\vspace{0.2cm}

2.) This direction can be found in \cite{CJR}. Let $C \subset S_1$ be
the exceptional curve of $\pi$ and $Z_1$ the support of
$\I_{\F_1}$. Then $C$ does not meet $Z_1$ and 
 \[F = f_1^{-1}(C) \simeq \PF_1\]
with minimal section $C_1$. We find $-K_{X_1}.C_1 = 0$ and the normal
bundle of $C_1$ in $X_1$ is of type $({-}1, {-}1)$. We may hence flop
$C_1$ and obtain $X_0^+$. The strict transform of $F$ is now a
contractible $\PN_2$, we denote the image of the blowdown by
$X_0$. Then $X_0$ is a $\PN_1$ bundle over the smooth almost Fano
surface $S_0$, hence $X_0 = \PN(\F_0)$ for some rank two vector bundle
$\F_0$. To prove $\F_0$ is of type \eqref{ext0} with $\I_{\F_0}$ as
claimed we chase a general member $H_1 \in |\sO_{X_1}(1)|$ as above.
\end{proof}

\vspace{0.2cm}

Back to our original situation. 
 First note that the blowup of $\PN_1\times\PN_1$ or $\PF_2$ in a
 general point is $\PN_2$ blown up in 2 points. To see this, start with
 $\PN_2$. Blowing up a point, we obtain $\PF_1$ with minimal section
 $C_0$. Choose some fiber $f$.
  \begin{enumerate}
    \item Blowing up a general point on $f$, we obtain three
      $({-}1)$--curves: the new exceptional divisor $E$, the strict
      transform $\hat{f}$ of $f$, and the minimal section $C_0$, which does not
      meet $E$. Blowing down $\hat{f}$ yields $\PN_1 \times \PN_1$.
    \item Blowing up the intersection point of $f$ and $C_0$ with
      exceptional divisor $E$, the strict transform $\hat{f}$ of $f$
      again is a $({-}1)$--curve, but now the strict transform
      $\hat{C}_0$ of $C_0$ is a  $({-}2)$--curve, not meeting
      $\hat{f}$. Blowing down $\hat{f}$ yields $\PF_2$. 
 \end{enumerate}
\vskip .2cm \noindent
Therefore we find a threefold diagram as
above (cf. \cite{CJR}):

\[\xymatrix{& X^+_1 \ar[r]^{\Bl} & X_1 \ar@{-->}[r]^{\text{\em flop}}
  \ar'[d][dd] &
  X^+_0 \ar[r]^{\Bl} & X_0 \ar[dd]\\ 
            X_2 \ar@{<-->}[rr]^{flop} \ar@{<-->}[ur]^{
\text{\em flop}}
            \ar[dd] && X^+_2 \ar[r]^{\Bl} & X \ar[dd]&\\
      &&{\PF_1} \ar'[r][rr] && \PN_2\\
            \Bl_p(\PF_2) \ar[urr] \ar[rrr] &&& \PF_2&}\]
(here formulated for $\PF_2$; we obtain exactly the same diagram for
$\PN_1 \times \PN_1$). Let $-K_{X_i} = 2H_i$. Then $H_0^3 = H_2^3+1$,
hence $H_0^3 \ge 2$. By the lemma, $X_0 = \PN(\F_0)$ for some rank two
vector bundle $\F_0$ of type \eqref{ext0}, where moreover the support
$\I_{\F_0}$ consists of the same points as the support of $\I_{\F}$ we
started with. It remains hence to classify all possible $X_0$ over
$\PN_2$. We find
 \begin{enumerate} 
   \item If $X_0$ is Fano, then either $\F_0 = T_{\PN_2}$, or $\F_0 =
     \sO(1) \oplus \sO(2)$ by Theorem~\ref{delpezzo} (the blowup of $\PN_3$ in a point, we have
     to normailze $\F_0$, such that $c_1(\F_0) = \sO(3)$). This gives $c_2(\F) = 3$, or $c_2(\F) = 2$,
     respectively.
  \item Assume $X_0$ is not Fano, the anticanonical map
    divisorial. Then $\F_0$ is one of the bundles in
    Theorem~\ref{div}, (2). We obtain $c_2(\F) = 6, 3, 0$ in cases
    (a), (b), (c), respectively.
 \item Assume $X_0$ is not Fano, the anticanonical map
    small. Then $\F_0$ is one of the bundles in
    Theorem~\ref{cb}, hence $4 \le c_2(\F) \le 7$.
\end{enumerate}

\

The existence of all cases is done in the following proposition.

\begin{proposition}  Let $X$ be a smooth almost
Fano threefold with $-K_X = 2H$ and $H^3 = d \ge 2$. Let $p \in X$ be
a general point and 
 \[Y = \Bl_p(X) \stackrel{\pi}{\lra} X.\]
Then $Y$ is again almost Fano.
\end{proposition} 

\begin{proof}
Let $E$ be the exceptional divisor of $\pi.$ Then $-K_Y = 2(\pi^*H -E) =: 2H'$. We get
$(H')^3 = H^3-1 > 0$. We claim $H'$ is nef.

By \cite{Shin}, $|H|$ is base point free, two general members $S_1,
S_2 \in |H|$ are hence smooth surfaces with $-K_{S_i} =
H|_{S_i}$. Their intersection is a smooth elliptic curve $C$. We may
assume $p \in C$. Blowing up $p$, we obtain $S'_i = \Bl_p(S_i)$ are the
strict transforms of $S_i$, hence $S'_i \in |H'|$. Moreover the
intersection $S'_1 \cap S'_2$ is isomorphic to $C$. To show $H'$ is
nef, it suffices to prove $H'|_{S'_1} = \sO_{S'_1}(C)$ is nef. Since
$C$ is irreducible, we have to prove the self intersection of $C$ in
$S'_1$ is non--negative. We have
 \[C\cdot_{S'_1}C = S'_1\cdot S'_1 \cdot S'_2 = (H')^3 = H^3-1 > 0.\]  
\end{proof}

\subsection*{The Vector Bundle Method} 
We consider a nef vector bundle $\sF$ over $S = \bP_1 \times \bP_1 $ (the case that $S$ is a blown-up quadric cone $\PF_2$ is very similar and
therefore omitted). Again we normalize $\sF$ such that
$$ \det \sF = -K_S.$$ 
Since $X = \bP(\sF)$ is assumed not to be Fano, $-K_{\bP(\sF)}$ is big and nef, but not ample, hence $\sF$ is not ample, but
$$c_1^2(\sF) > c_2(\sF).$$ 
We consider the anticanonical map $$\psi: X \to X'.$$ 
Then $H = \psi^*(H')$ and by [Sh89], $H'$ is always spanned unless $(H')^3 = 1,$ in which case $H'$ has a simple base point away from the
singularities of $X'.$ In all cases we find a section of $\sF$ vanishing in codimension at least $2$ (in most cases $F$ is even spanned). We thus have
an exact sequence
$$ 0 \to \sO_S \to \sF \to \sI_Z \otimes -K_S \to 0, \eqno (S)$$
where $Z$ 
is the zero locus of a general section of $\sF$ (so that the length $l(Z) = c_2 (\sF)$). 
The case $Z = \emptyset$ is very simple: here $(S)$ must split and thus $\sF = -K_S \oplus \sO_S,$ and this case of course really exists. 
The extension in case $l(Z) = 1$ does not exist (with $\sF$ locally free). This is easily seen by either restricting to ruling lines or by 
showing that $H^0(\sF(-1,-2)) \ne 0;$ so that $\sF(-1,-2)$ must have a section with at most finite zero locus; on the other hand 
$c_2(\sF(-1,-2)) = -1.$ \\  
Hence we are reduced to $l(Z) \geq 2.$ Since $K_S^2 = 8 $ and since $c_1^2(\sF) > c_2(\sF),$ we 
also have
$$ l(Z) \leq 7.$$ 
If $2 \leq l(Z) \leq 7$, then we study the restriction of $\sF$ to ruling lines $l_i$, i.e., to fibers of the projection $p_i: S \to B_i = \bP_1.$ 
Suppose that 
$$ \sF \vert l_2 = \sO(1) \oplus \sO(1)  $$
for all $l_2.$ Then $\sF(-1,0) \vert l_2 = \sO \oplus \sO, $ hence
$$ \sF(-1,0) = p_2^*(V) $$
with a vector bundle $V$ on $B_1.$ Hence $$ \sF = p_2^*(V) \otimes p_1^*(\sO(1)).$$
Identifying $l_1$ and $B_2,$ we see that $\sF \vert l_1 \simeq V,$ hence $V$ is nef
and either $V = \sO(1) \oplus \sO(1)$ or $V = \sO(2) \oplus \sO.$ The first alternative is impossible since then $V$ would be ample. 
So $\sF = \sO(1,2) \oplus \sO(1,0).$ 
\vskip .2cm \noindent
In a completely symmetric way, if $\sF \vert l_1$ for all $l_1,$ then $\sF = \sO(2,1) \oplus \sO(0,1).$
Thus we may assume that $\sF$ is not uniform on both ruling families. 
\vskip .2cm \noindent
Choose ``splitting lines'' $l_1^*$ and $l_2^*$ and $\tilde l_i$ be the exceptional sections in $\bP(\sF)$ sitting over $l_i^*.$ Then
$$ -K_X \cdot \tilde l_i = 0 $$
so that $\psi$ contracts $\tilde l_1$ and $\tilde l_2.$ This is only possible when two of the points of $Z$ are on $l_1^*$ or two of the points
ly on $l_2^*$. 
\vskip .2cm \noindent Conversely, take a finite set $Z \subset S$ with $2 \leq l(Z) \leq 7$ and with the distribution just described, the remaining 
points being in general position. Now ``Cayley-Bacharach for vector bundles `` (e.g. [GH78,p.731]) tells us that there is a vector bundle $\sF$
fitting into the exact sequence 
$$ 0 \to \sO_S \to \sF \to \sI_Z \otimes -K_S \to 0.$$ 
The special position of the points guarantees that $\sF$ is not ample. It remains to show that $\sF$ is nef. In fact, $\sF$ is spanned outside a finite
set, since the linear system $\vert \sI_Z \otimes -K_S \vert $ has no base components. 
\end{proof}

\begin{remark}
Note that $|H|$ is not base point free for $H^3 = 1$ by \cite{Shin}, but the
argument concerning blowups in the last paragraph of the proof also
applies, since two general members of $|H|$ still cut out an
irreducible curve. This means the blowup of any almost Fano threefold $X$ with $H^3=1$ in a
general point gives a smooth threefold $Y$ with $-K_Y$ nef, but not
big.
\end{remark}


\section{Almost del Pezzo manifolds in arbitrary dimension} 
\setcounter{lemma}{0}
In this section we consider an almost Fano manifold $X$ of dimension $n \geq 4$ admitting a line bundle $H$ such that 
$$ -K_X = (n-1)H.$$
So $X$ is an ``almost del Pezzo manifold''. We shall assume that $X$ is not Fano. Let $\phi: X \to Y$ be an extremal contraction; $\psi: X \to X'$
will again denote the anticanonical map to the singular del Pezzo variety $X'.$ Furthermore we have a line bundle $H'$ on $X'$ such that 
$-K_{X'} = (n-1)H'$. By [Fu90], we know that $H'$, hence $H$, is spanned if 
$ H^n \geq 2.$ We recall the classification of Gorenstein del Pezzo varieties [Fu90]:

\begin{proposition} \label{Fuclass} 
Let $X'$ be a Gorenstein del Pezzo $n-$fold $(n \geq 4$) of degree $d = (H')^n. $ Then $X'$ is one of the following. 
\begin{enumerate} 
\item $d = 1:$ \ $X'$ is a weighted hypersurface of degree $6$ in $\bP(3,2,1, \ldots, 1)$;
\item $d = 2:$ \ $X'$ is a weighted hypersurface of degree $4$ in $\bP(2,1, \ldots, 1);$ i.e., a double cover of $\bP_n$ branched along
a hypersurface of degree $4;$
\item $d = 3:$ \ $X' \subset \bP_{n+1}$ is a cubic;
\item $d = 4:$ \ $X' \subset \bP_{n+2}$ is a complete intersection of two quadrics;
\item $d \geq 5$ and $X'$ is a cone;
\item $d \geq 5$, \ $X'$ is not cone and $(n,d) = (4,6), (4,5), (5,5).$
\end{enumerate} 
\end{proposition}   

From results of Mori theory, e.g. [AW97,1.10,5.1] we obtain 

\begin{proposition} 
$\phi$ is one of the following.
\begin{enumerate}
\item a $\bP_{n-2}-$bundle over a smooth surface $Y$;
\item a quadric bundle over $\bP_1;$
\item the blow-up of a smooth point in the almost del Pezzo $n-$fold $Y$.
\end{enumerate}
\end{proposition} 

\subsection*{Case A: $\bP_{n-2}-$bundles}
We begin by treating case (1) and write
$$ X = \bP(\sF)$$
with $\sF$ a vector bundle of rank $n-1$ over $Y.$ As in the threefold case we can arrange $\det \sF = -K_Y$ and $-K_Y$ will be big and nef. 
If $d \geq 2,$ then $\sF$ is spanned by (5.1). Take $n-3$ general sections of $\sF$, then these sections give rise to an exact sequence
$$ 0 \to \sO_Y^{n-3} \to \sF \to \sF' \to 0 \eqno{(*)}$$
with a rank 2-bundle $\sF',$ cp. [OSS80,4.3.1]. 
$\sF'$ is spanned, in particular nef, and $\det \sF' = -K_Y.$ Thus $\bP(\sF')$ is an almost Fano threefold - possibly Fano - and therefore classified
by Theorem~\ref{div} and Theorem~\ref{dP}  resp. Theorem~\ref{delpezzo}. Notice also that $d = H^n = c_1^2(\sF) - c_2(\sF) = K_Y^2 - c_2(\sF) \leq 9.$ \\
Conversely, take a rank 2-bundle $\sF'$ which is nef but not ample; furthermore $c_1^2(\sF) > c_2(\sF),$ and $\det \sF' = -K_Y.$ Define $\sF$ by the exact
sequence
$$ 0 \to \sO_Y^{n-3} \to \sF \to \sF' \to 0,$$
e.g. $\sF = \sF' \oplus \sO_Y^{n-3}.$ Let $$X = \bP(\sF). $$ Then $X$ is almost del Pezzo. \\
Suppose now that $d = 1.$ Then by [Fu90,6.14], $\vert H' \vert$ has a simple base point $x_0$ which lies on the smooth part of $X'.$ Hence $\vert H \vert$
has just one simple base point $x_0.$ Let $y_0 = \pi(x_0)$, $\pi: X
\to Y$ the projection. In particular $\sF$ is generated outside $y_0$ and we obtain a sequence (*) on $Y \setminus y_0.$
This sequence is given by sections $s_1, \ldots, s_{n-3}$ which are linearly independent on $Y \setminus y_0,$ hence on $Y$. This (*) exists on all of
$Y$ and we can continue as before. 
We obtain:

\begin{theorem} \label{nbundle}
\begin{enumerate} 
\item Let $X_n$ be almost del Pezzo of the form $\bP(\sF)$ with a rank $(n-1)$-bundle $\sF$ over a smooth surface $Y$. After a suitable twist, $\det \sF = -K_Y.$ 
Then $Y$ is almost del Pezzo, and
$\sF$ fits into an exact sequence
$$ 0 \to \sO_Y^{n-3} \to \sF \to \sF' \to 0 $$
with a rank 2-bundle $\sF'$, and $\bP(\sF')$ is an almost del Pezzo threefold (classified in section 4). 
\item Let $Y$ be an almost del Pezzo surface and $\sF'$ a rank 2-bundle such that $\bP(\sF') $ is almost del Pezzo. Define $\sF$ as an extension
$$ 0 \to \sO_Y^{n-3} \to \sF \to \sF' \to 0,$$
e.g. $\sF = \sF' \oplus \sO_Y^{n-3}.$
Then $X = \bP(\sF)$ is almost del Pezzo, and not del Pezzo unless $X = \bP_2 \times \bP_2$ (with $\sF = \sO(1)^3 $ and $\sF' = T_{\bP_2}.$)
\end{enumerate}
\end{theorem} 

\begin{proof}
Only the last part of the theorem needs an explanation. Namely, assume $X = \bP(\sF)$ to be Fano. Then we use the classification to conclude that
$X = \bP_2 \times \bP_2.$ A priori it might happen that $\sF'$ is ample and that therefore also $\sF$ is ample. However the del Pezzo classification
shows that this can only happen when $\sF' = T_{\bP_2}.$
\end{proof} 
 
\subsection*{Case B: Quadric bundles} 
We now approach the second case, namely that 
$$\phi: X \to Y = \bP_1$$ is a quadric bundle. We introduce the rank $(n+1)-$vector bundle
$$ \sE = \phi_*(H)$$
so that $X \subset \bP(\sE)$. We also notice that $H = \sO_{\bP(\sE)}(1) \vert X.$ 

\begin{theorem} \label{cone}
If $X$ carries a quadric bundle structure, $X'$ cannot be a cone with
one exception: $\psi$ is small and the cone admits a second small resolution
$\tilde{X} = \PN(\sF)$, a $\PN_{n-2}$ bundle over $\PN_2$ as in
Theorem~\ref{nbundle} above. Using the same notation, $\sF$ is
determined by the rank two vector bundle $\sF'$ with $\sF'(2)$ as in (1) of Theorem~\ref{cb}.  
\end{theorem} 

\begin{proof} Suppose $X'$ is a cone. Then there is a birational map $f: \tilde X \to X'$, a del Pezzo variety $\tilde Z$ (which is not a cone)  carrying a vector bundle $V$
such that $$\tilde X = \bP(V)$$
with projection $p: \tilde X \to \tilde Z.$ Furthermore 
$\tilde Z$ is a general linear section by elements of $\vert H'
\vert$, hence Gorenstein with at most canonical singularities.   
Let $\tilde H = f^*(H');$ we normalize $V$ such that  
$$ \tilde H = \zeta_V := \sO_{\bP(V)}(1). $$ 
Let $r$ be the rank of $V$; then we can write
$$ -K_{\tilde X} = r \zeta_V + p^*(\det V^* \otimes -K_{\tilde Z}) = f^*((n-1)H') - \sum a_j E_j. \eqno (*)$$
Since $X'$ has only canonical singularities, so does $\tilde Z$, hence all $a_j \geq 0.$ 

\vskip .2cm \noindent {\bf (\ref{cone}.1)} 
First we assume that $f$ is small. Then the $E_j$ do not occur. 
Let $F$ be a general fiber of $p.$ 
Then, restricting (*) to $F,$ we obtain 
$$ r  =  n-1.$$ 
In that case $\dim \tilde Z = 2,$ hence $\tilde Z$ is a del Pezzo surface with canonical singularities. 
Let $h: \hat Z \to \tilde Z$ be the minimal desingularization so that $-K_{\hat Z} = h^*(-K_{\tilde Z}).$ 
Let $ \hat V = h^*(V)$ and set
$$ \hat X = \bP(h^*(V))= \tilde X \times_{\tilde Z}\hat Z$$
with projections $\hat h: \hat X \to \tilde X$ and $\hat{p}: \hat{X}
\to \hat{Z}$. We obtain
$$-K_{\hat X} = \hat{h}^*(-K_{\tilde X}) $$ 
so that $-K_{\hat X}$ is divisible by $n-1.$ Now $\hat Z $ admits a map $g: \hat Z \to \bP_1$ unless $\tilde Z = \bP_2$.
Then consider the  general fiber $G$ of $g  \circ \hat p$ and observe that
$-K_G$ is divisible by $n-1$, so that $G$ is a smooth quadric. On the other hand, $G$ admits a map to $\bP_1$, which yields a contradiction
since $n \geq 4.$ 

In the remaining case $\tilde{Z} \simeq \PN_2$ we note that $\tilde{X}
= \PN(V)$ is an almost del Pezzo manifold with $\rho(\tilde{X}) =
2$. These are classified in Theorem~\ref{nbundle}, i.e. there exists
an exact sequence
 \[0 \lra \sO_{\tilde{Z}}^{n-3} \lra V \lra V' \lra 0,\]
where $V'$ is a rank two vector bundle, such that $\tilde{X}_3 := \PN(V')$ is an
(almost) del Pezzo threefold and $c_1(V') = -K_{\tilde{Z}} = \sO_{\PN_2}(3)$. Moreover, $V = \sO^r \oplus M$ for some
rank $(n-1-r)$ bundle $M$ and $2r \le n$. 

If the threefold $\tilde{X}_3$ is not Fano, then the corresponding
anticanonical map $\tilde{X}_3 \to X'_3$ is small
and the flop $X_3$ admits a del Pezzo fibration. We conclude that $V'$
is one of the following list
 \begin{enumerate}
   \item $V' = \sF(2)$ with $\sF$ as in (1) of Theorem~\ref{cb},
   i.e. $c_1(\sF) = -1$ and $c_2(\sF) = 2$,
   \item $V' = T_{\PN_2}$,
   \item $V' = \sO(1) \oplus \sO(2)$.
 \end{enumerate}

\vskip .2cm \noindent 
(1) Assume $V' = \sF(2)$ with $\sF$ as in (1) of
Theorem~\ref{cb} and $r = n-3$, i.e. $V = V' \oplus \sO^{n-3}$. Then
$h^0(V(-1)) = h^0(\tilde{X}, H-p^*\sO_{\PN_2}(1)) = 2$ with base locus
excactly ${\rm exc}(f)$. Blowing up ${\rm exc}(f)$ and the contracting
the exceptional divisor the other direction first we obtain $X$
admitting a pencil. This is exactly the construction (1) in the proof of
Theorem~\ref{cb}.

\vskip .2cm \noindent 
(2,3) Concerning the other two cases, the quadric bundle structure of $X$ induces a
linear system $|L|$ on $\tilde X$ with exactly two sections. Let
 \[L = \alpha \eta_V - p^* \sO_{\PN_2}(\beta)\]
for some $\alpha, \beta \ge 0$. Then $h^0(S^{\alpha}V' \otimes
\sO(-\beta)) = 2$, which is impossible in the two remaining cases $V'=
T_{\PN_2}$ and $V' = \sO(1) \oplus \sO(2)$.

\vskip .2cm \noindent 
{\bf (\ref{cone}.2)} $f$ is divisorial. Then 
$$ V = \sO^r \oplus M $$
with a line bundle $M$ on $\tilde Z$ and $E = \bP(\sO^r)$ is the exceptional divisor of $f;$ the map $f$ is nothing than the blow-up of $X'$ along the
vertex $f(E) \simeq \bP_{r-1}.$  
By considering $Z_0 = \bP(M)$ (isomorphic to $\tilde Z$) and restricting $\tilde H = \zeta_V = f^*(H')$ to $Z_0$, it follows that 
$\tilde H \vert Z_0 = M$ which means $M = H' \vert \tilde Z.$\\
Now we consider a fiber of $\phi$, which is an $(n-1)-$dimensional quadric $Q_{n-1}$ and take its $\psi-$image $Q'_{n-1}$ which is isomorphic to $Q_{n-1}$
and which contains $f(E).$ Let $\tilde Q_{n-1}$ be the strict transform of $Q'_{n-1}$ in $\tilde X$; then 
$$\tilde Q_{n-1} \to Q'_{n-1}$$
is nothing than the blow-up of the smooth quadric $Q'_{n-1}$ along the linear subspace $f(E) = \bP_{r-1}.$ Now the blow-up of $\bP_n$ along $\bP_{r-1}$
is Fano with second projection to $\bP_{n-r}$. Hence also $\tilde Q_{n-1}$, a divisor in the blow-up of $\bP_n,$ has a surjective map to $\bP_{n-r}$ 
(with connected fibers; just the second projection of the Fano manifold $\tilde Q_{n-1}$).
Thus
$$ \tilde Z = \bP_{n-r}, $$
e.g. because $\overline{NE}(\tilde X)$ is a 2-dimensional cone. Since $\tilde Z$ is a linear section in $X'$, we obtain $-K_{\tilde{Z}} =
(n-1-r)H'|_{\tilde{Z}}$ by adjunction. Then 
 \[n-r+1 = (n-r-1)a\] 
with $H'|_{\tilde{Z}} = \sO(a)$, hence $n-1-r = 1$ and $a = 3$, or $n-1-r =
2$ and $a = 1$. Assume $n-1-r = 1$ and $a = 3$. Then $\tilde{Z} = \PN_2$ and $V = \sO_{\PN_2}^{n-2} \oplus \sO_{\PN_2}(3)$. 
Now we use again the fact that $\tilde X$ carries some line bundle with exactly two
sections. This is impossible. \\ 
If $n-1-r = 2$ and $a = 1$, then $\tilde{Z} = \PN_3$ and $V = \sO_{\PN_3}^{n-3} \oplus
\sO_{\PN_3}(1)$, so that $X' = \PN_n$, which is absurd. 
\end{proof}

We are considering next the special cases (6) in Proposition~\ref{Fuclass}.

\begin{proposition} The case $(n,d) = (4,6)$ does not occur.
\end{proposition}

\begin{proof} By \cite{Fu90} $X'$ is obtained in the following way. We consider the vector bundle
$$V = \sO_{\bP_2} \oplus \sO_{\bP_2}^{\oplus 3} $$
with projection $p: \bP(V) \to \bP_2$ 
and let $\zeta = \sO_{\bP(V)}(1).$
Then $\vert \zeta \vert$ defines a morphism
$$ \pi: \bP(V) \to \bP_8 $$ 
contracting the divisor 
$$D = \bP(\sO_{\bP_2}^{\oplus 3}) \simeq \bP_2 \times \bP_2$$
to $\bP_2.$ Notice that $D \in \vert \zeta - p^*(\sO(2)) \vert.$ 
Now  $$ \tilde X \in \vert \zeta + p^*(\sO(1)) \vert $$
is a general member, and $X' = \pi(\tilde X)$ with induced map
$\tilde{\pi}: \tilde{X} \to X'$. Obviously $\tilde D = D \cap \tilde X$ is a divisor of type $(1,1)$ 
in $D \simeq \bP_2 \times \bP_2,$ hence
$$\tilde D \simeq \bP(T_{\bP_2})$$
and therefore $\tilde{\pi}(\tilde D) \simeq \bP_2$ in $X'.$ Thus $\tilde{\pi}$ is divisorial and ${\rm Sing}(X') \simeq \bP_2.$ 
The adjunction formula shows that $-K_{\tilde X} = 3 \zeta,$ hence $\tilde{\pi}$ is crepant, so that $X'$ is $\bQ-$factorial with canonical
non-terminal singularities. Hence $\psi: X \to X'$ defined by $\vert H \vert$ cannot be small, otherwise $X'$ would not be $\bQ-$factorial. 
So $\psi$ contracts an irreducible divisor $E$ to $\bP_2$ (recall that $X$ carries a quadric bundles structure so that $\rho(X) = 2$). 
Notice also that $\psi$ has connected fibers (otherwise consider the Stein factorisation yielding a covering $X'' \to X'$ with $X''$ singular
del Pezzo which cannot exist by Fujita's classification). 
Now let $Q_3 $ be a general fiber of $\phi.$ Then $\psi \vert Q_3$ is finite and has degree $1$. Let $Q_3' = \psi(Q_3) \subset  X'.$ 
>From the adjunction formula we see that
$$H \vert Q_3 = \sO_{Q_3}(1), $$
hence $H^0(X,H) \to H^0(Q_3,H \vert Q_3)$ must be surjective. This means that $Q_3'$ is a quadric in $\bP_4.$
Since on the other hand $\psi \vert Q_3$
is generically an isomorphism and since $E \cap Q_3$ is a divisor on $Q_3,$ it follows that $\psi(E) \subset Q_3'.$ But an irreducible
quadric in $\bP_4$ cannot contain a $\bP_2.$ 
\end{proof}

\begin{theorem} \label{case55}
The case $(n,d) = (5,5)$ occurs: there exists an almost del Pezzo $5-$fold $X$ of degree $5$ admitting a quadric bundle
structure over $\bP_1.$ Moreover $X$ is not Fano.
\end{theorem} 

\begin{proof} We first give the description of $X'$ as in \cite{Fu90}, (9.9.5), (9.14) and (9.12.si). Let $\sE = \sO_{\bP_1}(1)^{\oplus 3}$ and
set $$M = \bP(\sE) \simeq \bP_1 \times \bP_2 $$
with projection $p: M \to \bP_1$ and ``tautological'' line bundle $\sO_M(1).$
Let 
$$ W = \bP(\sO_M(1) \oplus \sO_M^3) $$
with tautological bundle $\zeta$ and projection $f: W \to M$.   
Let 
$$ \tilde X \in \vert \zeta + f^*(\sO_M(1) - F) \vert $$
be general, where $F$ is a fiber of $p.$ Since $ \zeta + f^*(\sO_M(1) - F)$ is clearly spanned, $\tilde X$ is smooth. 
Then $\zeta $ defines a map $\pi: W \to \bP_8$ which is birational onto its image with exceptional divisor  
$$D \in \vert \zeta - f^*(\sO_M(1)) \vert.$$ 
The variety $X'$ is just the $\pi-$image of $\tilde X:$
$$ X' = \pi(\tilde X) \subset \bP_8.$$ 
The adjunction formula gives
$$ K_{\tilde X} = -3 \zeta_{\tilde X} - f^*(\sO_M(1)) \vert \tilde X =
- 4 \zeta_{\tilde X} + D \vert \tilde X. \eqno (*)$$
Since $\zeta = \pi^*(\sO_{\bP_8}(1)) $ and since $K_{X'} = \sO_{\bP_8}(-4) \vert X',$ it follows
$$ K_{\tilde X} = \tilde \pi^*(K_{X'}) + \tilde D, $$
where $\tilde \pi = \pi \vert \tilde X$ and $\tilde D = D \vert \tilde X.$ 
Hence $X'$ has only terminal singularities. \\
We have a closer look to the $5-$fold $D.$ Since $D = \bP(\sO_M^3), $ we have $D \simeq \bP_1 \times \bP_2 \times \bP_2$, and, in order to keep
track of the projective plane, we write more specifically
$$D = \bP_1 \times \bP_2^a \times \bP_2^b,$$ 
where $M = \bP_1 \times \bP_2^a.$ Hence $\pi(D) = \bP_2^b$ and $f(D) = \bP_1 \times \bP_2^a.$ 
Now $\tilde X \vert D$ is a divisor of type $(0,1,1)$ hence 
$$ \tilde D = \bP_1 \times \bP(T_{\bP_2}). $$ It follows that the singular locus of $X'$ is $\bP_2^b.$ \\
If $X$ exists, then $\rho (X') = 1. $ On the other hand, $\rho(\tilde X) = 3.$ Hence we try to factorize $\tilde \pi$ and to obtain
$X$ as intermediate variety. \\
The line bundle $\zeta + f^*(F) $ is spanned and big, but clearly not ample. Let 
$$ g: W\to W''$$
be the associated birational morphism which is clearly divisorial. 
If $W' = \pi(W),$ then $\pi: W \to W'$ factorizes as
$$ W {\buildrel {g} \over \to } W'' {\buildrel {h} \over \to } W'.$$ 
The exceptional divisor of $g$ is still $D$ and $g(D) = \bP_1 \times \bP_2^b$, whereas $h(\bP_1  \times \bP_2^b) =  \bP_2^b.$ 
Intersecting with $\tilde X,$ we obtain a birational map $\tilde g: \tilde X \to g(\tilde X)$ such that $\tilde g$ contracts 
$\tilde D = \bP_1 \times \bP(T_{\bP_2}) $ to $\bP_1  \times \bP_2^b$ and then $h$ projects to $\bP_2^b.$   
We set 
$$ X = g(\tilde X). $$
We need to show that $X$ is smooth, del Pezzo, and admits a quadric
bundle structure. The smoothness is seen as follows. Take a line $l$
in a $\tilde g$-fiber. 
Then 
$$ \tilde D \cdot l = D \cdot l = \zeta \cdot (- f^*(\sO_M(1))) = -1.$$
Hence $$\tilde D \vert \tilde g^{-1}(x) = \sO(-1) $$
and Nakano's theorem says that $X$ is smooth. The divisibility of $K_X$ comes from (*). 
Finally the quadric bundle is induced from $\tilde X \to \bP_1,$ since
$\tilde g$ only contracts curves in 
fibers of $\tilde X \to \bP_1$ (observe that $ \zeta + f^*(F)$ is ample on all curves projecting onto $\bP_1$).  
\end{proof}

\begin{theorem} Suppose $(n,d) = (4,5).$ Then $X$ is a hyperplane section of a del Pezzo manifold of type $(5,5).$ 
\end{theorem} 

\begin{proof} Applying again \cite{Fu90}, we are either in case 9.14(7). Here all the computations of Proposition~\ref{case55} work in the same way and
it is clear that we obtain a hyperplane section of a $5-$fold of type $(5,5).$ \\
Or we are in case 9.14(6) of \cite{Fu90}. We are going to rule out this case. The description of $X'$ is very similar to that one in (\ref{case55}); we are going
to use the same notations.  
Here we consider the Hirzebruch surface
$$ p: M = \bP(\sO(2) \oplus \sO(1)) \to \bP_1$$ 
and set
$$ W = \bP(\sO_M(1) \oplus \sO_M^3)$$
with birational map $\pi: W \to W'$ provided by $\zeta.$ The exceptional divisor is
$$ D \in \vert \zeta - f^*(\sO_M(1)) \vert,$$ 
so that $D = M \times \bP_2 = \PF_1 \times \bP_2.$ 
We take
$$\tilde X \in \vert \zeta + f^*(\sO_M(1) - F) \vert $$
general, $X' = \pi(\tilde X).$ Then 
$$ K_{\tilde X} = -3 \zeta$$
and thus $\tilde \pi: \tilde X \to X'$ is crepant and divisorial. In particular $X'$ has canonical non-terminal singularities. 
We consider the exceptional divisor
$$ \tilde D = D \cap \tilde X$$
of $\tilde \pi.$ Inside $D = \PF_1 \times \bP_2$ it is of type
$(C_0+F,1),$ where $C_0$ is the $(-1)-$curve in $\PF_1.$ Then $\tilde
\pi(\tilde D) = \bP_2$
is the singular locus of $X'.$ \\
Suppose now that $X$ exists. Then $\rho(X') = 1$ and we must have a factorization
$$ \tilde X {\buildrel {g} \over \to } X'' {\buildrel {h} \over \to} X'.$$ 
Now $\tilde D$ is a $\bP_1$-bundle over $\PF_1,$ so that $\rho(\tilde
D) = 3.$ Hence $\tilde \pi: \tilde D \to \bP_2$ must have singular
fibers, and $g \vert \tilde D$ will 
contract components of singular fibers, so that $g(\tilde D)$ is a
$\bP_1-$bundle over $\bP_2.$ In other words, $g \vert \tilde D$ is just the restricition of
the blow-down map $D = \PF_1 \times \bP_2 \to \bP_2 \times \bP_2.$ The conclusion is that $g$ is small, while $h$ is divisorial. 
Hence $X''$ is not $\bQ-$factorial and so does $X'.$ \\
In summary $X'$ is neither terminal nor $\bQ-$factorial. But our potential $X$ has $\rho(X) = 2.$ Hence either $\psi: X \to X'$ is small - then $X'$ 
would be terminal. Or $\psi$ is divisorial - then $X'$ would be $\bQ-$factorial. This leads to the contradiction we are looking for, and $X$ cannot
exist. 

\end{proof}

We summarize the results in Case B:

\begin{theorem} Let $X$ be an almost del Pezzo manifold which is a quadric bundle over $\bP_1.$ Let $\psi: X \to X'$ be the anticanonical map.
\begin{enumerate} 
\item $X'$ is never a cone with the following exception: $\psi$ is small and $X'$ admits a small resolution $\tilde X \to X$ such that 
$\tilde X$ is a $\bP_{n-2}-$bundle over $\bP_2$ as decribed in (5.3). The associated rank 2-bundle $\sF'(2)$ is given in (3.6)(1).
\item $\dim X = 5, \ H^5 = 5$ and $X'$ is a del Pezzo 5-fold of degree $5,$ i.e. a singular hyperplane section of $G(1,4).$
\item $\dim X = 4, \ H^5 = 4$ and $X'$ is a hyperplane section of (2). 
\end{enumerate}
All cases really exist. 
\end{theorem}

\subsection*{Case C: Blow-ups}
The case that $\phi: X \to Y$ is the blow-up of a smooth point is settled by 

\begin{theorem} (1) Let $X_n$ be almost del Pezzo, $n \geq 3,$ and $\phi: X \to Y$ be the blow-up of a point in the manifold $Y$. 
Then $Y$ is almost del Pezzo. If $-K_Y = (n-1)H,$ then $H^n > 1.$  \\
(2) Conversely, let $Y_n$ be del Pezzo and $\phi: X \to Y$ be the blow-up of a {\it general} point $p \in Y.$ Write 
$-K_Y = (n-1)H$ and assume $H^n > 1.$ 
Then $X$ is almost del Pezzo. 
\end{theorem} 
 
\begin{proof} (1) This is completely obvious using 
$$ -K_X = \phi^*(-K_Y) - (n-1)E, $$
where $E$ is the exceptional divisor. \\
(2) In the other direction we proceed by induction on $n,$ the case $n = 3$ being settled by Proposition 4.5.
So let $n \geq 4.$  
By [Fu90], (3.5) and (4.16), applied to the anticanonical model $Y'$
of $Y$, the line bundle $H$ is spanned, since $H^n > 1.$ Let $S \in \vert H \vert$ be a smooth member. 
Then $S$ is an almost del Pezzo $(n-1)-$fold, and we may assume that $p \in S$, $p$ being general. Let $\hat S$ be the strict transform of $S$ in
$X$, the blow-up of $S$ at $p.$ So by induction $\hat S$ is again a
del Pezzo $(n-1)-$fold. Write $-K_X = (n-1)\hat H,$ so that $\hat S
\in \vert \hat H \vert $.
Since $-K_{\hat S} = (n-2)\hat H \vert \hat S$, the line bundle $\hat
H \vert \hat S$ is nef. Thus $\hat H$ itself is nef and so does $-K_X.$ Since 
$(\hat H)^n = H^n - 1 > 0,$ the manifold $X$ is del Pezzo. 
\end{proof} 



\begin{thebibliography}{OSS80}
\bibitem[AW97]{AW97} M. Andreatta, J. Wisniewski: A view on
  contraction of higher dimensional varieties. Proc. Symp. Pure
  Math. {\bf 62.1}, 153-183 (1997).
\bibitem[BS95]{BS95} M. Beltrametti, A. Sommese: The adjunction theory
  of complex projective varieties. De Gruyter exp. in math. {\bf 16} (1995).
\bibitem[Bo01]{Borisov} A. Borisov: Boundedness of Fano threefolds
  with log-terminal singularities of given index. J. Math. Sci., Tokyo
  {\bf 8}, 329-342 (2001).
\bibitem[CJR06]{CJR} C. Casagrande, P. Jahnke, I. Radloff: On the
  Picard number of almost Fano threefolds with pseudo-index $>1$. math.AG/0603090.
\bibitem[Cu88]{Cu} S. Cutkosky: Elementary Contractions of Gorenstein Threefolds. Math. Ann. {\bf 280}, 521-525 (1988).
\bibitem[Fu80]{Fu80} T. Fujita: On the structure of polarized manifolds
  with total deficiency one, 1,2 and 3. J. Math. Soc. Japan {\bf 32},
  709-725 (1980); {\bf 33}, 415-434 (1981); {\bf 36}, 75-89 (1984). 
\bibitem[Fu90]{Fu90} T. Fujita: Classification theories of polarized
  varieties. London Math. Lect. Notes {\bf 155}, Cam. Univ. Press (1990).
\bibitem[GH78]{GH78} Griffiths,Ph., Harris,J.: Principles of Algebraic Geometry. Wiley (1978)
\bibitem[Hu79]{Hulek} K. Hulek: Stable Rank-2 Vector bundles on
  $\PN_2$ with $c_1$ odd. Math. Ann. {\bf 242}, 241-266 (1979)
\bibitem[I78]{Isk78} V.A. Iskovskikh: Fano $3$-folds I, II. Math. USSR, Izv. {\bf 11}, 485-527 (1977); {\bf 12}, 469-506 (1978).
\bibitem[I80]{Isk80} V.A. Iskovskikh: Anticanonical models of
  three-dimensional algebraic varieties. J. Soviet Math. {\bf 13},
  745-814 (1980). 
\bibitem[IP99]{AG5} V.A. Iskovskikh, Yu.G. Prokhorov: Algebraic
  Geometry V: Fano varieties. Springer 1999.
\bibitem[JPR05]{JPR} P. Jahnke, T. Peternell, I. Radloff: Threefolds
  with big and nef anticanonical bundles I. Math. Ann. {\bf 333},
  569-631 (2005). 
\bibitem[JR06]{smoothing} P. Jahnke, I. Radloff: Terminal Fano
  threefolds and their smoothings. math.AG/0601769.
\bibitem[Ka88]{Kaw} Y. Kawamata: Crepant blowing-up of
  $3$-dimensional canonical singularities and its application to
  degeneration of surfaces. Ann. Math. {\bf 127}, 93-163 (1988).
\bibitem[Ko89]{Kollar} J. Koll\'ar: Flops. Nagoya Math. J. {\bf 113}, 15-36 (1989).
\bibitem[KM98]{KoMo} J. Koll\'ar, S. Mori: Birational Geometry of
  Algebraic Varieties. Camb. Univ. Press 1998.
\bibitem[LP79]{LeP79} J. Le Potier: Fibr\'es stables de rang $2$ sur
  $\PN_2(\KC)$. Math. Ann. {\bf 241}, 217-256 (1979)
\bibitem[MK02]{McK} J. McKernan: Boundedness of log terminal Fano pairs
  of bounded index. math.AG/0205214.
\bibitem[Mi01]{Mi} T. Minagawa: Deformations of weak Fano $3$-folds
  with only terminal singularities. Osaka J. Math. {\bf 38}, 533-540 (2001).
\bibitem[Mo82]{Mori} S. Mori: Threefolds whose canonical bundles are not numerically effective. Ann. Math.
 {\bf 116}, 133-176 (1982).
\bibitem[Na97]{Namikawa} Y. Namikawa: Smoothing Fano $3$-folds. J. Alg. Geom. {\bf 6}, 307-324 (1997).
\bibitem[OSS80]{OSS} C. Okonek, M. Scneider, H. Spindler: Vector
  bundles on complex projective spaces. Progress in math. {\bf 3},
  Birkh\"auser 1980.
\bibitem[Pr05]{Prokh} Yu.G. Prokhorov: The degree of Fano threefolds
  with canonical Gorenstein singularities. Mat.\ Sb\ {\bf 196}, 81-122
  (2005), in Russian. English translation: Sb.\ Math.\ {\bf 196},
  77-114 (2005).
\bibitem[Re83]{Reid} M. Reid: Projective morphisms according to
  Kawamata. Unpublished manuscript (1983).
\bibitem[Sh89]{Shin} K.-H. Shin: $3$-Dimensional Fano varieties with
    canonical singularities. Tokyo J. Math. {\bf 12}, 375-385 (1989).
\end{thebibliography}
\end{document}